\documentclass[11pt, leqno]{amsart}
\usepackage{amssymb}
\usepackage{amsmath,amscd}
\usepackage[initials,nobysame,alphabetic]{amsrefs}
\usepackage{amsmath, amssymb}
\usepackage{amsfonts}
\usepackage{mathrsfs}
\usepackage{mathpazo}
\usepackage[arrow,matrix,curve,cmtip,ps]{xy}

\usepackage{amsthm}

\usepackage{float}
\usepackage{hyperref}
\usepackage{graphics}
\usepackage{graphicx}
\usepackage{verbatim}
\usepackage{multirow}
\usepackage{tikz}
\usepackage{enumerate}
\usepackage{eufrak}
\allowdisplaybreaks

\newtheorem{theorem}{Theorem}[section]
\newtheorem{lemma}[theorem]{Lemma}
\newtheorem{proposition}[theorem]{Proposition}

\newtheorem*{theorem*}{Theorem}
\newtheorem{remark}[theorem]{Remark}
\newtheorem{definition}[theorem]{Definition}
\newtheorem{example}[theorem]{Example}

\numberwithin{equation}{section}

\newcommand{\U}{\mathcal{U}}

\newcommand{\Prob}{\mathbb{\Prob}}

\newcommand{\mytilde}{\raise.17ex\hbox{$\scriptstyle\mathtt{\sim}$}}

\usepackage{amsmath}
\usepackage{amsfonts}
\usepackage{enumerate}
\usepackage{amssymb}
\usepackage{graphicx}%
\providecommand{\U}[1]{\protect\rule{.1in}{.1in}}

\begin{document}

\title{Mean-field FBSDE and optimal control}
\author{Nacira Agram and Salah Eddine Choutri}

\address{Department of Mathematics, Linnaeus University, Vaxjo, Sweden. Email: nacira.agram@lnu.se.}
\address{Department of Mathematics, KTH Royal Institute of Technology,
100 44, Stockholm, Sweden. Email: choutri@kth.se.}

\thanks{This research was carried out with support of the Norwegian
Research Council, within the research project Challenges in Stochastic
Control, Information and Applications (STOCONINF), project number
250768/F20.}
\date{This version \today}



\begin{abstract}
We study optimal control for mean-field forward backward stochastic
differential equations with payoff functionals of mean-field type. Sufficient and necessary optimality conditions in terms of a stochastic maximum principle are derived. As
an illustration, we solve an optimal portfolio with mean-field risk minimization problem.
\end{abstract}

\maketitle
\tableofcontents

\section{Introduction}

\noindent Stochastic differential equation (SDE) of mean-field type (a.k.a McKean-Vlasov equation) is an SDE whose coefficients depend on the marginal law of the solution (state) as well as the solution itself, i.e.
\small 
\begin{equation} \label{mean-field}
\left\{
\begin{array}
[c]{ll}
dX(t) & =b(t,X(t),
\mathcal{L}(X(t)))dt+\sigma(t,X(t),
\mathcal{L}(X(t)))dB(t),\\
X(0) & =x_{0}.
\end{array}
\right.
\end{equation} \normalsize
where $\mathcal{L}(X(t))$ is the law of $X(t)$, which is obtained as a limit of a sequence of empirical distribution functions representing the states. \\
The SDE \eqref{mean-field} can be viewed as the limit of a system of particles with mean-field interaction  \small
\[
\left\{
\begin{array}
[c]{ll}%
dX^{i,n}(t) & =b(t,X^{i,n}(t),%

\frac{1}{n}%
{\textstyle\sum_{j=1}^{n}}
\delta_{X^{j,n}(t)}%

)dt +\sigma(t,X^{j,n} (t),%

\frac{1}{n}%
{\textstyle\sum_{j=1}^{n}}
\delta_{X^{j,n}(t)}

)dB^{i}(t),\\
X(0) & =x_{0},
\end{array}
\right.
\]
\normalsize
 when the size of the system $n$ tends to infinity. \\
Optimal control of mean-field SDEs was first studied  by
Andersson and Djehiche \cite{andersson2011}, and Buckdahn, Djehiche and Li \cite{buckdahngeneral}, where the mean-field coupling is represented by an expected value of the state. The authors established a suitably modified stochastic maximum principle which involves mean-field backward SDEs.\\
\noindent Extension to the case where the marginal law of the state process is the mean-field coupling was studied by Carmona and Delarue \cite{carmona1}. The authors use the Wassertein metric
space for measures and the lifting technique introduced by
Lions \cite{Lions} to differentiate a function of a measure. In \cite{AO1}, Agram and \O ksendal have introduced a Sobolev space of random measures in which, the Fr\'{e}chet derivative with respect to the measure can be used directly. This approach is used in the present paper.

\noindent The purpose of our work is to derive necessary and sufficient optimality conditions in terms of a stochastic maximum principle for a set $\hat{u}$ of admissible controls which minimize a cost functional of the form
\[%
\begin{array}
[c]{ll}%
J(u) & =E[h(X(T),M(T))+\phi(Y(0),N(0))\\
& +{%
{\textstyle\int_{0}^{T}}
}f(t,X(t),Y(t),Z(t),M(t),N(t),u(t))dt],
\end{array}
\]
with respect to admissible controls $u$, for some functions $f,\ h,\ \phi$, under dynamics governed by mean-field forward backward stochastic
differential equations (MF-FBSDE). More specifically, we consider the coupled system
\small
\begin{equation*}
\left\{
\begin{array}
[c]{lll}%
dX(t) & = & b(t,X(t),M(t),u(t))dt+\sigma(t,X(t),M(t),u(t))dB(t),\text{ \ }%
t\in\left[  0,T\right]  ,\\
X(0) & = & x_{0},\text{ \ }%
\end{array}
\right.  \label{eq:X}%
\end{equation*}

\begin{equation*}
\left\{
\begin{array}
[c]{lll}%
dY(t) & = & -g(t,X(t),Y(t),Z(t)),M(t),N(t),u(t))dt+Z(t)dB(t),\text{ \ }%
t\in\left[  0,T\right]  ,\\
Y(T) & = & X(T),
\end{array}
\right.  \label{eq:Y}%
\end{equation*} \normalsize
for some functions $b,\ \sigma$ and a Brownian motion $B(t).$ $M(t)$ and $N(t)$ denote the marginal laws of X and Y respectively.

\noindent Existence of a fully-coupled MF-FBSDE is
studied by Carmona and Delarue under Lipschitz assumption on the coefficients but no
uniqueness result was proven. Bensoussan et al \cite{bens} prove existence and
uniqueness of a fully coupled MF-FBSDE by assuming Lipschitz and monotonicity conditions. Recently, Djehiche and Hamadene in \cite{Boualem2019} prove the same results but under weak monotonicity assumptions and without the non-degeneracy condition on the forward equation.

\noindent In the next section, we give some mathematical background. Next, we study stochastic optimal control of
MF-FBSDE where sufficient and necessary optimality conditions are derived. In the
last section, we construct a discounted dynamic risk measure by means of
MF-BSDE and then we solve an associated risk minimization problem.

\section{Generalities}

\noindent Let $B=B(t),t\in\lbrack0,T]$ be a one-dimensional Brownian motion
defined in a complete filtered probability space $(\Omega,\mathcal{F}%
,\mathbb{F},P).$ The filtration $\mathbb{F=}\left\{  \mathcal{F}%
_{t}\right\}  _{t\geq0}$ is assumed to be the $P$-augmented
filtration generated by $B.$\newline

\begin{definition}
\-

\begin{itemize}
\item Let $\mathcal{M}$ be the space of random measures $\mu$ on $\mathbb{R}$
equipped with the norm
\begin{equation}%
\begin{array}
[c]{lll}%
\left\Vert \mu\right\Vert _{\mathcal{M}}^{2} & := & E[%
{\textstyle\int_{\mathbb{R}}}
|\hat{\mu}(y)|^{2}e^{-y^{2}}dy]\text{,}%
\end{array}
\label{f_nrm}%
\end{equation}
where $\hat{\mu}$ is the Fourier transform of the measure $\mu$, i.e.,%
\[%
\begin{array}
[c]{lll}%
\hat{\mu}(y) & := & {%
{\textstyle\int_{\mathbb{R}}}
}e^{ixy}d\mu(x);\quad y\in\mathbb{R}.
\end{array}
\]
We endow $\mathcal{M}$ with the inner product $\langle\mu,\eta\rangle:=%
{\textstyle\int_{\mathbb{R}}}
|\hat{\mu}(y)-\hat{\eta}(y)|^{2}e^{-y^{2}}dy,$ $\mu,\eta,y\in\mathbb{R}$,
$\hat{\mu}$ and $\hat{\eta}$ are the Fourier transform of the measures $\mu$
and $\eta$. Then $(\mathcal{M},||\cdot||)$ is a pre-Hilbert space.

\item We denote by $\mathcal{M}_{0}$ the set of all deterministic elements of
$\mathcal{M}$.
\end{itemize}
\end{definition}

We give some examples:

\begin{example}
[Measures]Let us give some examples of measures in $\mathcal{M}_{0}$ and
$\mathcal{M}$:\-

\begin{enumerate}
\item Suppose that $\mu=\delta_{x_{0}}$, the unit point mass at $x_{0}%
\in\mathbb{R}$. Then $\delta_{x_{0}}\in\mathcal{M}_{0}$ and
\[%
\begin{array}
[c]{lll}%
\hat{\mu}(y) & ={%
{\textstyle\int_{\mathbb{R}}}
}e^{ixy}d\mu(x) & =e^{ix_{0}y},
\end{array}
\]
and hence
\[%
\begin{array}
[c]{lll}%
\left\Vert \mu\right\Vert _{\mathcal{M}_{0}}^{2} & =%
{\textstyle\int_{\mathbb{R}}}
|e^{ix_{0}y}|^{2}e^{-y^{2}}dy & <\infty\text{.}%
\end{array}
\]

\item Suppose $d\mu(x)=f(x)dx$, where $f\in L^{1}(\mathbb{R})$. Then $\mu
\in\mathcal{M}_{0}$ and by Riemann-Lebesque lemma, $\hat{\mu}(y)\in
C_{0}(\mathbb{R})$, i.e. $\hat{\mu}$ is continuous and $\hat{\mu
}(y)\rightarrow0$ when $|y|\rightarrow\infty$. In particular, $|\hat{\mu}|$ is
bounded on $\mathbb{R}$ and hence%
\[%
\begin{array}
[c]{lll}%
\left\Vert \mu\right\Vert _{\mathcal{M}_{0}}^{2} & =%
{\textstyle\int_{\mathbb{R}}}
|\hat{\mu}(y)|^{2}e^{-y^{2}}dy & <\infty\text{.}%
\end{array}
\]

\item Suppose that $\mu$ is any finite positive measure on $\mathbb{R}$. Then
$\mu\in\mathcal{M}_{0}$ and
\[%
\begin{array}
[c]{lll}%
|\hat{\mu}(y)| & \leq%
{\textstyle\int_{\mathbb{R}}}
d\mu(y)=\mu(\mathbb{R}) & <\infty\text{, for all }y\text{,}%
\end{array}
\]
and hence%
\[%
\begin{array}
[c]{lll}%
\left\Vert \mu\right\Vert _{\mathcal{M}_{0}}^{2} & =%
{\textstyle\int_{\mathbb{R}}}
|\hat{\mu}(y)|^{2}e^{-y^{2}}dy & <\infty\text{.}%
\end{array}
\]

\item Next, suppose $x_{0}=x_{0}(\omega)$ is random. Then $\delta
_{x_{0}(\omega)}$ is a random measure in $\mathcal{M}$. Similarly, if
$f(x)=f(x,\omega)$ is random, then $d\mu(x,\omega)=f(x,\omega)dx$ is a random
measure in $\mathcal{M}$.\newline
\end{enumerate}
\end{example}

\noindent We denote by $U$ a nonempty convex subset of $\mathbb{R}$ and we
denote by $\mathcal{U}_{\mathbb{G}}$ the set of $U$-valued $\mathbb{G}%
$-progressively measurable processes where $\mathbb{G}:=\{\mathcal{G}%
_{t}\}_{t\geq0}$ with $\mathcal{G}_{t}\subseteq\mathcal{F}_{t}$ for all
$t\geq0$; we consider them as the admissible control processes.\newline We
will also use the following spaces:

\begin{itemize}
\item $\mathcal{S}^{2}$ is the set of ${\mathbb{R}}$-valued $\mathbb{F}%
$-adapted c\`{a}dl\`{a}g processes $X=X(t),t\in\lbrack0,T],$ such that
\[
{\Vert X\Vert}_{\mathcal{S}^{2}}^{2}:=E[\sup_{t\in\lbrack
0,T]}|X(t)|^{2}]~<~\infty\;,
\]

\item $\mathbb{L}^{2}$ is the set of ${\mathbb{R}}$-valued $\mathbb{F}%
$-adapted processes $Q=Q(t),t\in\lbrack0,T],$ such that
\[
\Vert Q\Vert_{\mathbb{L}^{2}}^{2}:=E[%
{\textstyle\int_{0}^{T}}
|Q(t)|^{2}dt]<~\infty\;.
\]

\item $\mathcal{K}$ denotes the set of absolutely continuous functions
$m:[0,T]\rightarrow\mathcal{M}_{0}.$

\item $\mathbb{K}$ is the set of bounded linear functionals $K:\mathcal{M}%
_{0}\rightarrow\mathbb{R}$ equipped with the operator norm
\[
||K||_{\mathbb{K}}:=\sup_{m\in\mathcal{M}_{0},||m||_{\mathcal{M}_{0}}\leq
1}|K(m)|.
\]

\item $\mathcal{S}_{\mathbb{K}}^{2}$ is the set of $\mathbb{F}$-adapted
stochastic processes $p:[0,T]\times\Omega\mapsto\mathbb{K},$ such that
\[
||p||_{\mathcal{S}_{\mathbb{K}}}^{2}:=E[\sup_{t\in\lbrack
0,T]}||p(t)||_{\mathbb{K}}^{2}]<\infty.
\]

\item $\mathbb{L}_{\mathbb{K}}^{2}$ is the set of $\mathbb{F}$-adapted
stochastic processes $q:[0,T]\times\Omega\mapsto\mathbb{K},$ such that
\[
||q||_{\mathbb{L}_{\mathbb{K}}^{2}}^{2}:=E[%
{\textstyle\int_{0}^{T}}
||q(t)||_{\mathbb{K}}^{2}dt]<\infty.
\]

\end{itemize}

\noindent We recall now the notion of differentiability which will be used in
the sequel.\newline Let $\mathcal{X},\mathcal{Y}$ be two Banach spaces with
norms $\Vert\cdot\Vert_{\mathcal{X}},\Vert\cdot\Vert_{\mathcal{Y}}$,
respectively, and let $F:\mathcal{X}\rightarrow\mathcal{Y}$.

\begin{itemize}
\item We say that $F$ has a directional derivative (or Gateaux derivative) at
$v\in\mathcal{X}$ in the direction $w\in\mathcal{X}$ if
\[
D_{w}F(v):=\lim_{\varepsilon\rightarrow0}\frac{1}{\varepsilon}(F(v+\varepsilon
w)-F(v))
\]
exists in $\mathcal{Y}$.

\item We say that $F$ is Fr\'{e}chet differentiable at $v\in\mathcal{X}$ if
there exists a continuous linear map $A:\mathcal{X}\rightarrow\mathcal{Y}$
such that
\[
\lim_{\substack{h\rightarrow0\\h\in\mathcal{X}}}\frac{1}{\Vert h\Vert
_{\mathcal{X}}}\Vert F(v+h)-F(v)-A(h)\Vert_{\mathcal{Y}}=0,
\]
where $A(h)=\langle A,h\rangle$ is the action of the liner operator $A$ on
$h$. In this case we call $A$ the \textit{gradient} (or Fr\'{e}chet
derivative) of $F$ at $v$ and we write
\[
A=\nabla_{v}F.
\]

\item If $F$ is Fr\'{e}chet differentiable at $v$ with Fr\'{e}chet derivative
$\nabla_{v}F$, then $F$ has a directional derivative in all directions
$w\in\mathcal{X}$ and
\[
D_{w}F(v)=\nabla_{v}F(w)=\langle\nabla_{v}F,w\rangle.
\]

\end{itemize}

\noindent In particular, note that if $F$ is a linear operator, then
$\nabla_{v}F=F$ for all $v$.\newline

\section{Optimal control problem}

\noindent Here we denote by $M(t):=\mathcal{L}(X(t))$ the law of $X(t)$ at
time $t$ and by $N(t):=\mathcal{L}(Y(t))$ the law of $Y(t)$ at time $t$. We
assume that our system is gouverned by a coupled system of MF-FBSDE as
follows:\newline The MF-SDE $X^{u}(t)=X(t)$ is given by \small
\begin{equation}
\left\{
\begin{array}
[c]{lll}%
dX(t) & = & b(t,X(t),M(t),u(t))dt+\sigma(t,X(t),M(t),u(t))dB(t),\text{ \ }%
t\in\left[  0,T\right]  ,\\
X(0) & = & x_{0},\text{ \ }%
\end{array}
\right.  \label{eq:X}%
\end{equation} \normalsize
for functions $\sigma,b:\Omega\times\left[  0,T\right]  \times%
\mathbb{R}
\times\mathcal{M}_{0}\times U\rightarrow%
\mathbb{R}
$ which are supposed to be $\mathcal{F}_{t}$-measurable and the initial value
$x_{0}\in%
\mathbb{R}
$.\newline The couple MF-BSDE $(Y^{u}(t),Z^{u}(t))=(Y(t),Z(t))$ satisfies \small
\begin{equation}
\left\{
\begin{array}
[c]{lll}%
dY(t) & = & -g(t,X(t),Y(t),Z(t)),M(t),N(t),u(t))dt+Z(t)dB(t),\text{ \ }%
t\in\left[  0,T\right]  ,\\
Y(T) & = & X(T),
\end{array}
\right.  \label{eq:Y}%
\end{equation} \normalsize
where $g:\Omega\times\left[  0,T\right]  \times%
\mathbb{R}
^{3}\times\mathcal{M}_{0}^{2}\times U\rightarrow%
\mathbb{R}
$ is $\mathbb{F}$-adapted.\newline

\noindent It is obvious from the definition of the norm (\ref{f_nrm}) that%
\[%
\begin{array}
[c]{lll}%
||\mathcal{L}(X^{(1)})-\mathcal{L}(X^{(2)})||_{\mathcal{M}_{0}}^{2} & \leq &
\sqrt{\pi}\mathbb{E}[(X^{(1)}-X^{(2)})^{2}]\text{,}%
\end{array}
\]
where $X^{(1)}$ and $X^{(2)}$ are random variables that follow the
distributions $\mathcal{L}(X^{(1)})$ and $\mathcal{L}(X^{(2)})$ respectively.

\noindent Assume that ($C$ is a constant that may change
from line to line)

\begin{itemize}
\item[(A1)] there exists $C>0$, such that

\item for all $t\in\left[  0,T\right]  $, for all fixed $u\in U,$
$x,x^{\prime}\in%
\mathbb{R}
,m,m^{\prime}\in\mathcal{M}_{0}$%
\begin{align*}
\left\vert \sigma\left(  t,x,m,u\right)  -\sigma\left(  t,x^{\prime}%
,m^{\prime},u\right)  \right\vert +\left\vert b\left(  t,x,m,u\right)
-b\left(  t,x^{\prime},m^{\prime},u\right)  \right\vert \\
\leq C\left(\left\vert x-x^{\prime}\right\vert+||m-m^{\prime}||_{\mathcal{M}_{0}}\right).
\end{align*}

\item for all $t\in\left[  0,T\right]  $,\ for all fixed $u\in U,$%
\[
\left\vert \sigma\left(  t,0,\delta_{0},u\right)  \right\vert +\left\vert
b\left(  t,0,\delta_{0},u\right)  \right\vert \leq C\text{,}%
\]
where $\delta_{0}$ is the distribution law of zero, i.e., the Dirac measure
with mass at zero.

\item[(A2)] there exists $C>0$, such that, for all fixed $u\in U$ and all
knowing $X(t)\in\mathcal{S}^{2}$ of equation (\ref{eq:X}) and $M(t):= \mathcal{L}%
(X(t))\in\mathcal{M}_{0}$, we have

\item for all $t\in\left[  0,T\right]  ,$ $y,y^{\prime},z,z^{\prime}\in%
\mathbb{R}
,n,n^{\prime}\in\mathcal{M}_{0}$

\begin{align*}
\left\vert g\left(  t,x,y,z,m,n,u\right)  -g\left(  t,x,y^{\prime},z^{\prime
},m,n^{\prime},u\right)  \right\vert \\ 
\leq C\left(  \left\vert y-y^{\prime}\right\vert +\left\vert z-z^{\prime}\right\vert +||n-n^{\prime}%
||_{\mathcal{M}_{0}}\right)  .
\end{align*}

\item for all $t\in\left[  0,T\right]  $,
\[
\left\vert g\left(  t,x,0,0,m,\delta_{0},u\right)  \right\vert \leq C\text{.}%
\]

\end{itemize}

\begin{proposition}
Under Assumptions (A1) and (A2), the MF-FBSDE (\ref{eq:X})-(\ref{eq:Y}) admits
a unique solution $(X,Y,Z)\in\mathcal{S}^{2}\times\mathcal{S}^{2}%
\times\mathbb{L}^{2}.$
\end{proposition}

\noindent Since the system is partially-coupled i.e., the forward equation does
not depend on the solution of the backward one, we can solve the system
separately as follows: we first find a solution $X(t)$ of the MF-SDE (\ref{eq:X})
and then we plug it into the backward equation \eqref{eq:Y}, then we solve it. 

\noindent Our aim is to maximize the performance functional of the form%
\[%
\begin{array}
[c]{ll}%
J(u) & =E[h(X(T),M(T))+\phi(Y(0),N(0))\\
& +{%
{\textstyle\int_{0}^{T}}
}f(t,X(t),Y(t),Z(t),M(t),N(t),u(t))dt],
\end{array}
\]
over all admissible controls, for functions $f:\Omega\times\left[  0,T\right]
\times%
\mathbb{R}
^{3}\times\mathcal{M}_{0}^{2}\times U\rightarrow%
\mathbb{R}
$, $h:\Omega\times%
\mathbb{R}
\times\mathcal{M}_{0}\rightarrow%
\mathbb{R}
$ and $\phi:\Omega\times%
\mathbb{R}
\times\mathcal{M}_{0}\rightarrow%
\mathbb{R}
.$\newline Now, we can define the Hamiltonian%
\[
H:\Omega\times\left[  0,T\right]  \times%
\mathbb{R}
^{3}\times\mathcal{M}_{0}^{2}\times U\times%
\mathbb{R}
^{2}\times\mathcal{K}\times%
\mathbb{R}
\times\mathcal{K}\rightarrow%
\mathbb{R}
\]
by%

\begin{align}
H(t,x,y,z,m,n,u,p^{0},q^{0},p^{1},\lambda^{0},\lambda^{1})  &
=f(t,x,y,z,n,u)+p^{0}b(t,x,m,u)\nonumber\\
&  +q^{0}\sigma(t,x,m,u)+\lambda^{0}g(t,x,n,u)\nonumber\\
&  +\langle p^{1},m^{\prime}\rangle+\langle\lambda^{1},n^{\prime}\rangle.
\label{H}%
\end{align}

\begin{remark}
For ease of notation we drop the dependence of all variables except for the
time $t,\forall\Phi\in\{\sigma,f,H,h,g,\phi\},$ we write $\Phi(t), \forall t.$
\ Moreover, we will use
\begin{align*}
\hat{\Phi}(t):  &  =\Phi(t,\hat{X}(t),\hat{Y}(t),\hat{Z}(t),\hat
{M}(t),\hat{N}(t),\hat{u}(t))\\
\check{\Phi}(t) :  &  =\Phi(t,\hat{X}(t),\hat{Y}(t),\hat{Z}(t),\hat{M}(t),\hat{N}(t),u(t))].
\end{align*}

\end{remark}

\noindent For $u\in\mathcal{U}$ with corresponding solution $X^{u}=X$, define,
whenever solutions exist, $p^{\hat{u}}=p=(p^{0},p^{1})$ and $q^{\hat{u}%
}=q=(q^{0},q^{1})$ and $\lambda^{\hat{u}}=\lambda=(\lambda^{0},\lambda^{1})$
by the adjoint equations:\newline The BSDE for the unknown processes
$(p^{0},q^{0})\in\mathcal{S}^{2}\times\mathbb{L}^{2}$
\begin{equation}
\left\{
\begin{array}
[c]{lll}%
dp^{0}(t) & = & -\partial_{x}H(t)dt+q^{0}(t)dB(t),t\in\left[  0,T\right]  ,\\
p^{0}(T) & = & \partial_{x}h(T)+\lambda^{0}(T).
\end{array}
\right.  \label{eq:P0}%
\end{equation}
The MF-BSDE for the unknown processes $(p^{1},q^{1})\in\mathcal{S}%
_{\mathbb{K}}^{2}\times\mathbb{L}_{\mathbb{K}}^{2}$
\begin{equation}
\left\{
\begin{array}
[c]{lll}%
dp^{1}(t) & = & -\nabla_{m}H(t)dt+q^{1}(t)dB(t),t\in\left[  0,T\right]  ,\\
p^{1}(T) & = & \nabla_{m}h(T),
\end{array}
\right.  \label{eq:P1}%
\end{equation}
The forward SDE%
\begin{equation}
\left\{
\begin{array}
[c]{lll}%
d\lambda^{0}(t) & = & -\partial_{y}H(t)dt+\partial_{z}H(t)dB(t),t\in\left[
0,T\right]  ,\\
\lambda^{0}(0) & = & \partial_{y}\phi(0),
\end{array}
\right.  \label{eq:lambda 0}%
\end{equation}
and
\begin{equation}
\left\{
\begin{array}
[c]{lll}%
d\lambda^{1}(t) & = & -\nabla_{n}H(t)dt,t\in\left[  0,T\right]  ,\\
\lambda^{1}(0) & = & \nabla_{n}\phi(0).
\end{array}
\right.  \label{eq:lambda 1}%
\end{equation}

\noindent Before stating and proving sufficient and necessary conditions of
optimality, we need the following result, which is Lemma 2.3 in Agram and \O ksendal
\cite{AO1}:

\begin{lemma}
Suppose that $X(t)$ is an It\^{o} process of the form
\[%
\begin{cases}
dX(t)=\theta(t)dt+\gamma(t)dB(t),\quad t\in\lbrack0,T],\\
X(0)=x_{0}\in\mathbb{R},
\end{cases}
\]
where $\theta,\gamma$ are adapted processes.\newline Then the map $t\mapsto
M(t):[0,T]\rightarrow\mathcal{M}_{0}$ is absolutely continuous.
\end{lemma}

\noindent It follows that $t\mapsto M(t)$ is differentiable for $t$-a.e. We
will in the following use the notation
\[
M^{\prime}(t)=\tfrac{d}{dt}M(t).
\]

\subsection{Sufficient optimality conditions}

\begin{theorem}
Suppose that $\hat{u}\in\mathcal{U}_{\mathbb{G}}$ with corresponding solutions
$\hat{X}(t),(\hat{Y}(t),\hat{Z}(t)), \newline(p^{0}(t),q^{0}(t)),(p^{1}%
(t),q^{1}(t)),\lambda^{0}(t),\lambda^{1}(t)$ to equations \eqref{eq:X},
\eqref{eq:Y}, \eqref{eq:P0},\eqref{eq:P1}, \eqref{eq:lambda 0} and
\eqref{eq:lambda 1} respectively. Suppose that

\begin{itemize}
\item $x,m\mapsto h(k,m)$ ,

\item $y,n\mapsto\phi(y,n)$ ,

\item $x,y,z,m,n,u\mapsto H(\cdot,x,y,z,m,n,u)$ ,
\end{itemize}

are concave functions $P$-a.s for each $t \in[0,T]. $ \ Moreover,%

\[
E[\hat{H}(t)|\mathcal{G}_{t}]=\max_{u\in U}E[\check{H}(t)|\mathcal{G}_{t}],
\]
$P$-a.s for all t $\in\lbrack0,T].$ \noindent Then $\hat{u}$ is an optimal control.
\end{theorem}

\noindent{Proof}\quad We show that $J(u)-J(\hat{u})\leq0,$ for an arbitrary
$u$ and a fixed \newline optimal $\hat{u} \in\mathcal{U}_{\mathbb{G}}. $ \newline

\noindent We introduce first the following notation $\forall\Phi\in
\{\sigma,f,H,h,g,\phi,M,N,M^{\prime},N^{\prime}\}$ and $\forall t,$%
\[
\delta\Phi(t)=\check{\Phi}(t)-\hat{\Phi}(t),
\]
and
\[
\delta M^{\prime}(t)=\delta(\frac{d}{dt}M(t))=\frac{d}{dt} (\delta M(t)).
\]
From the definition of the Hamiltonian (\ref{H}), we have%

\begin{align*}
\delta f(t) &  =\delta H(t)-\delta b(t)p^{0}(t)-\delta\sigma(t)q^{0}(t)\\
&  -\langle p^{1}(t),M^{\prime}(t)\rangle-\langle\lambda^{1}(t),N^{\prime
}(t)\rangle,
\end{align*}
and
\begin{align}
J(u)-J(\hat{u}) &  =E[%
{\textstyle\int_{0}^{T}}
\{\delta H(t)-\delta b(t)p^{0}(t)-\delta\sigma(t)q^{0}(t)-\langle
p^{1}(t),M^{\prime}(t)\rangle\label{J-J}\\
&  -\langle\lambda^{1}(t),N^{\prime}(t)\rangle\}dt+\delta h(T)+\delta
\Phi(0)].\nonumber
\end{align}
We use the concavity of $h$ and $\phi$ as well as the boundary values of equations
\eqref{eq:P0}, \eqref{eq:P1}, \eqref{eq:lambda 0} and \eqref{eq:lambda 1}%

\begin{align}
\delta h(T)+\delta\phi(0) &  \leq\partial_{x}h(T)\ \delta X(T)+\langle
\nabla_{m}h(T),\delta M(T)\rangle  \nonumber  \\ \nonumber
&  +\partial_{x}\phi(0)\delta Y(0)+\langle\nabla_{n}\phi(0),\delta
N(0)\rangle   \\ \label{eq: concavity Phi and h}
&  =p^{0}(T)\delta X(T)-\lambda^0(T) \delta X(T)+\langle p^{1}(T),\delta M(T)\rangle\\ \nonumber
&  +\lambda^{0}(0)\delta Y(0)+\langle \lambda^{1}(0) ,\delta N(0) \rangle.
\end{align}
Applying It$\hat{o}$ formula to $p^{0}(t)\delta X(t),p^{0}(t)\delta
X(t),\lambda^{0}(t)\delta Y(t)$ and $\lambda^{1}(t)\delta Y(t),$ yields the
following duality relations:%

\begin{align} \label{Dual-sufficient 1}
E[p^{0}(T)\delta X(T))]=E[%
{\textstyle\int_{0}^{T}}
p^{0}(t)\delta b(t)dt-%
{\textstyle\int_{0}^{T}}
\delta X(t)\partial_{x}H(t)dt+%
{\textstyle\int_{0}^{T}}
Z(t)\delta\sigma(t)dt\ ],
\end{align}

\begin{align} \label{Dual-sufficient 2}
E[\langle p^{1}(T),\delta M(T)\rangle]=E[%
{\textstyle\int_{0}^{T}}
\langle p^{1}(t),\delta M^{\prime}(t)\rangle dt-%
{\textstyle\int_{0}^{T}}
\langle\nabla_{m} \hat{H} (t),\delta M(t)\rangle dt],
\end{align}

\begin{align} \label{Dual-sufficient 3}
E[\lambda^0(T) \delta Y(T) ]-E[\lambda^{0}(0)\delta Y(0)]   = &-E[%
{\textstyle\int_{0}^{T}}
\lambda^{0}(t)\delta g(t)dt]+E[%
{\textstyle\int_{0}^{T}}
\delta Y(t)\partial_{y} \hat{H}(t)dt]     \\ \nonumber
&  +E[%
{\textstyle\int_{0}^{T}}
Z(t)\partial_{z}\hat{H}(t) dt],
\end{align}

\begin{equation} \label{Dual-sufficient 4}
E[\lambda^1(T) \delta N(T) ]-E[\langle \lambda^{1}(0),\delta N(0)\rangle]=E[%
{\textstyle\int_{0}^{T}}
\langle\lambda^{1}(t),\delta N^{\prime}(t)\rangle+\langle\nabla
_{n}\hat{H}(t),\delta N(t)\rangle dt].
\end{equation}

\noindent By the concavity of $H,$ we obtain%

\begin{align} \label{suffcient-concavity hamiltonian}
\delta H(t) &  \leq\partial_{x}\hat{H}(t)\delta X(t)+\partial_{y}\hat
{H}(t)\delta Y(t)+\partial_{z}\hat{H}(t)\delta Z(t)\\ \nonumber 
&+\langle\nabla_{m}\hat{H}(t),\delta
M(t)\rangle+\langle\nabla_{m}\hat{H}(t),\delta N(t)\rangle+\partial_{u}\hat{H}(t)\delta u(t).
\end{align}

\noindent Finally, by substituting the derived duality relations \eqref{Dual-sufficient 1},\eqref{Dual-sufficient 2}, \eqref{Dual-sufficient 3} and \eqref{Dual-sufficient 4} in \eqref{J-J} and
using the estimates \eqref{eq: concavity Phi and h}, \eqref{suffcient-concavity hamiltonian}, we obtain%

\[
J(u)-J(\hat{u})\leq E[%
{\textstyle\int_{0}^{T}}
\partial_{u}\hat{H}(t)\delta u(t)]
\]
Using the tower property and the fact that $u(t)$ is $\mathbb{G}$-adapted the
desired result follows
\begin{align*}
J(u)-J(\hat{u}) &  \leq E[%
{\textstyle\int_{0}^{T}}
E[\partial_{u}\hat{H}(t)|\mathcal{G}_{t}]\ \delta u(t)dt]
  \leq0,
\end{align*}
and thus, $\hat{u}$ is optimal.

$\qquad\qquad\qquad\qquad\qquad\qquad\qquad\qquad\qquad\qquad\qquad
\qquad\qquad\qquad\qquad\square$ \newline

\subsection{Necessary optimality conditions}

\noindent Given an arbitrary but fixed control $u\in\mathcal{U}_{\mathbb{G}}$,
we define
\begin{equation} \label{perturbed u}
u^{\rho}:=\hat{u}+\rho u, \ \ \ \  \rho \in [0,1]  .
\end{equation}
Note that, the convexity of $U$ and $\mathcal{U}_{\mathbb{G}}$ guarantees that
$u^{\rho}\in\mathcal{U}_{\mathbb{G}},\rho\in\left[  0,1\right]  $. We denote
by $X^{\rho}:=X^{u^{\rho}}$ and by $\hat{X}:=X^{\hat{u}},$ the solution
processes corresponding to $u^{\rho}$ and $\hat{u},$\ respectively.\newline

\noindent For each $t_{0}\in\left[  0,T\right]  $ and all bounded $\mathcal{G}_{t_{0}}%
$-measurable random variables $\alpha,$ the process
\[%
\begin{array}
[c]{lll}%
u\left(  t\right)   & = & \alpha\mathbf{1}_{\left(  t_{0},T\right]  }(t),
\end{array}
\]
belongs to $\mathcal{U}_{\mathbb{G}}$.\medskip

\noindent In general, if $K^{\hat{u}}(t)$ is a process depending on $\hat{u}$, we define
the operator $D$ on $K$ by
\begin{equation}
DK^{\hat{u}}(t):=D^{u}K^{\hat{u}}(t)=\tfrac{d}{d\rho}K^{\hat{u}+\rho
u}(t)|_{\rho=0},
\end{equation}
whenever the derivative exists.\newline Define the following derivative
processes
\begin{align*}
DX^\rho (t) &:= \tfrac{d}{d\rho}X^{\hat{u}+\rho u}(t)|_{\rho=0}=\mathcal{X}^{\rho}(t) ,\\
DY^\rho (t) &:= \tfrac{d}{d\rho}Y^{\hat{u}+\rho u}(t)|_{\rho=0}=\mathcal{Y}^{\rho}(t) ,\\
DZ^\rho (t)&:= \tfrac{d}{d\rho}Z^{\hat{u}+\rho u}(t)|_{\rho=0}=\mathcal{Z}^{\rho}(t) ,\\
DN^{\rho}(t) &:= \tfrac{d}{d\rho}N^{\hat{u}+\rho u}(t)|_{\rho=0},\\
DM^{\rho}(t) &:= \tfrac{d}{d\rho}M^{\hat{u}+\rho u}(t)|_{\rho=0},\\
DN^{\rho^{\prime}}(t)&:= \tfrac{d}{d\rho}\tfrac{d}{dt}M^{\hat{u}+\rho u}(t)|_{\rho=0},\\
DM^{\rho^{\prime}}(t) &:= \tfrac{d}{d\rho}\tfrac{d}{dt}M^{\hat{u}+\rho u}(t)|_{\rho=0},
\end{align*}
such that
%

\begin{equation}
\left\{
\begin{array}
[c]{ll}%
d\mathcal{X}^{\rho}(t) & =\{\partial_{x}b(t)\mathcal{X}^{\rho}(t)+\nabla
_{m}b(t)DM^{\rho}(t)+\partial_{u}b(t)u(t)\}dt\\
& +\{\partial_{x}\sigma(t)\mathcal{X}^{\rho}(t)+\langle\nabla_{m}%
\sigma(t),DM^{\rho}(t)\rangle+\partial_{u}\sigma(t) u(t)\}dB(t),t\in\left[
0,T\right]  ,\\
\mathcal{X}^{\rho}(0) & =0,
\end{array}
\right.  \label{eq:perturbed X}%
\end{equation}
and
\begin{equation}
\left\{
\begin{array}
[c]{ll}%
d\mathcal{Y}^{\rho}(t) & =-\{\partial_{x}g(t)\mathcal{X}^{\rho}(t)+\partial
_{y}g(t)\mathcal{Y}^{\rho}(t)+\partial_{z}g(t)\mathcal{Z}^{\rho}%
(t)+\langle\nabla_{m}g(t),DM^{\rho}(t)\rangle\\
& +\langle\nabla_{n}g(t),DN^{\rho}(t)\rangle+\partial_{u}g(t) u
(t)\}dt+\mathcal{Z}^{\rho}(t)dB(t),t\in\left[  0,T\right]  ,\\
\mathcal{Y}^{\rho}(T) & =0.
\end{array}
\right.  \label{eq:perturbed Y}%
\end{equation}

\noindent Moreover, we assume that all the partial derivatives of $\Phi\in
\{\sigma,f,H,h,g,\phi\}$ are bounded.

\begin{theorem}
\label{Th: necessary} Let $\hat{u}\in\mathcal{U}_{\mathbb{G}}$ be the optimal
control and $\mathcal{X}^{\rho}(t),(\mathcal{Y}^{\rho}(t),\mathcal{Z}^{\rho}(t)),\newline (p^{0}(t),
q^{0}(t)),(p^{1}(t),q^{1}(t)),\lambda^{0}(t),\lambda
^{1}(t)$ be the corresponding solutions to the \newline equations
\eqref{eq:perturbed X},\eqref{eq:perturbed Y}, \eqref{eq:P0},\eqref{eq:P1},
\eqref{eq:lambda 0},\eqref{eq:lambda 1}. Then, the following statements are
\newline equivalent

\begin{enumerate}
[(i)]

\item $\frac{d}{d\rho}J(\hat{u}+\rho u)|_{\rho=0}=0$ for all bounded $\beta
\in\mathcal{U}_{\mathbb{G}}.$

\item $E[\frac{\partial}{\partial u}\hat{H}(t)|\mathcal{G}_{t}]=0$ for all
$t\in\lbrack0,T].$
\end{enumerate}
\end{theorem}

\noindent{Proof}\quad We first prove theorem \ref{Th: necessary} by assuming (i)
and aiming to show (ii)%

\begin{align*}
0 &  =\tfrac{d}{d\rho}J(u+\rho u)|_{\rho=0}\\
&  =E[%
{\textstyle\int_{0}^{T}}
\tfrac{d}{d\rho}\ f(t)|_{\rho=0}dt+p^{0}(T)\mathcal{X}^{\rho}(T)+\langle
p^{1}(T),DM^{\rho}(T)\rangle+\lambda^{0}(0)\mathcal{Y}^{\rho}(0)\\
&  +\langle\lambda^{1}(0),DN^{\rho}(0)\rangle]
\end{align*}
\{ we substitute $f(t)$ from equation $\eqref{H}$ \}
\begin{align*}
&  =E[%
{\textstyle\int_{0}^{T}}
\tfrac{d}{d\rho}\{H^\rho(t)-p^{0}(t)b^\rho(t)-q^{0}(t)\sigma^{\rho}(t)-\lambda
^{0}(t)g^\rho(t)-\langle p^{1}(t),M^{\rho^{\prime}}(t)\rangle\\
&  -\langle\lambda^{1}(t),N^{\rho^{\prime}}(t)\rangle\}|_{\rho=0}%
dt+p^{0}(T)\mathcal{X}^{\rho}(T)-\lambda^0(T) \mathcal{X}^{\rho}(T)+\langle p^{1}(T),DM^{\rho}(T)\rangle\\
& +\lambda^{0}(0)\mathcal{Y}^{\rho}(0)
+\langle\lambda^{1}(0),DN^{\rho}(0)\rangle],
\end{align*}
by using the chain rule, we obtain
\begin{align*}
\tfrac{d}{d\rho}H^\rho(t)|_{\rho=0} &  =\partial_{x}H(t)\mathcal{X}^{\rho
}(t)+\partial_{y}H(t)\mathcal{Y}^{\rho}(t)+\partial_{z}H(t)\mathcal{Z}^{\rho
}(t)+\langle\nabla_{m}H(t),DM^{\rho}(t)\rangle\\&+\langle\nabla_{n}H(t),DN^{\rho
}(t)\rangle+\partial_{u}H(t)u(t),
\end{align*}

\[
\tfrac{d}{d\rho}p^{0}(t)b^\rho(t)|_{\rho=0}=p^{0}(t)\partial_{x}b(t)\mathcal{X}%
^{\rho}(t)+p^{0}(t)\langle\nabla
_{m}b(t),DM^{\rho}(t)\rangle+p^{0}(t)\partial_{u}b(t)u(t),
\]

\[
\tfrac{d}{d\rho}q^{0}(t)\sigma^\rho(t)|_{\rho=0}=q^{0}(t)\partial_{x}%
\sigma(t)\mathcal{X}^{\rho}(t)+q^{0}(t)\langle\nabla_{m}\sigma(t),DM^{\rho
}(t)\rangle+q^{0}(t)\partial_{u}\sigma(t)u(t),
\]

\begin{align*}
\tfrac{d}{d\rho}\lambda^{0}(t)g^\rho(t)|_{\rho=0} &  =\lambda^{0}(t)\partial
_{x}g(t)\mathcal{X}^{\rho}(t)+\lambda^{0}(t)\partial_{y}g(t)\mathcal{Y}^{\rho
}(t)+\lambda^{0}(t)\partial_{z}g(t)\mathcal{Z}^{\rho}(t)\\
&  +\lambda^{0}(t)\langle\nabla_{m}g(t),DM^{\rho}(t)\rangle+\lambda
^{0}(t)\langle\nabla_{n}g(t),DN^{\rho}(t)\rangle\\
&  +\lambda^{0}(t)\partial_{u}g(t) u(t),
\end{align*}

\[
\tfrac{d}{d\rho}\langle p^{1}(t),M^{\rho^{\prime}}(t)\rangle|_{\rho=0}=\langle
p^{1}(t),DM^{\rho^{\prime}}(t)\rangle,
\]
\noindent and
\[
\tfrac{d}{d\rho}\langle\lambda^{1}(t),N^{\rho^{\prime}}(t)\rangle|_{\rho
=0}=\langle\lambda^{1}(t),DN^{\rho^{\prime}}(t)\rangle.
\]
We apply It$\hat{o}$ formula to $p^{0}(t)\mathcal{X}^{\rho}(t),\langle
p^{1}(t),DM^{\rho}(t)\rangle,\lambda^{0}(t)\mathcal{Y}^{\rho}(t)$ and \\
$\langle\lambda^{1}(t),DN^{\rho}(t)\rangle$ then we take the expectation, we
obtain the following important duality relations:

\begin{align*}
E[p^{0}(T)\mathcal{X}^{\rho}(T)]  &  =E[%
{\textstyle\int_{0}^{T}}
\{p^{0}(t)\partial_{x}b(t)\mathcal{X}^{\rho}(t)
+p^{0}(t)\langle\nabla_{m}b(t),DM^{\rho}(t)\rangle\\
&+p^{0}(t)\partial_{u}%
b(t)u(t)-\partial_{x}H(t)\mathcal{X}^{\rho}(t)+q^{0}(t)\partial_{x}\sigma
(t)\mathcal{X}^{\rho}(t)\\
& +q^{0}(t)\langle\nabla_{m}\sigma(t),DM^{\rho
}(t)\rangle
  +q^{0}(t)\partial_{u}\sigma(t) u(t)\}dt],
\end{align*}

\begin{align*}
E[\langle p^{1}(T),DM^{\rho}(T)\rangle]  &  =E[%
{\textstyle\int_{0}^{T}}
\langle p^{1}(t),DM^{\rho^{\prime}}(t) \rangle- \langle \nabla_{m}H(t), DM^{\rho}(t) \rangle dt],
\end{align*}

\begin{align*}
E[\lambda^0(T)\mathcal{Y}^\rho(T)]-E[\lambda^{0}(0)\mathcal{Y}^{\rho}(0)]  &  =E[%
{\textstyle\int_{0}^{T}}
\{ -\lambda^{0}(t)\partial_{x}g(t)\mathcal{X}^{\rho}(t)\\
&-\lambda^{0}(t)\partial_{y}g(t)\mathcal{Y}^{\rho}(t) 
-\lambda^{0}(t)\partial_{z}g(t)\mathcal{Z}^{\rho}(t) \\
&-\lambda^{0}(t)\langle\nabla_{m}g(t),DM^{\rho}(t)\rangle\\
&-\lambda^{0}(t)\langle\nabla_{n}g(t),DN^{\rho}(t)\rangle\\ 
&-\lambda^{0}(t)\partial_{u}g(t)u(t) +\partial_{y}H(t)\mathcal{Y}^{\rho}\\
&+\partial_{z}H(t)\mathcal{Z}^{\rho}(t)\}dt],
\end{align*}

\begin{align*}
E[\langle \lambda^1(T) , DN^\rho(T) \rangle]-E[\langle\lambda^{1}(0),DN^{\rho}(0)\rangle] &=E[{%
{\textstyle\int_{0}^{T}}
}\{\langle\lambda^{1}(t),DN^{\rho^{\prime}}(t)\rangle \\
&+\langle\nabla_{n}H(t),DN^{\rho}(t)\rangle\}dt].
\end{align*}
\newline By substituting the derived duality relations and the partial
derivatives of $f(t)$ the desired result follows. This proof can be reversed
to prove $(ii)\Rightarrow(i).$ We omit the details. $ \qquad\qquad\ \qquad\qquad\ \ \ \ \ \qquad \qquad\qquad\ \qquad \qquad\qquad\square$

\section{Mean-field discounted risk measure}

\noindent In this section we are interested in a particular class of MF-BSDE of the
following form 
\begin{equation}
\left\{
\begin{array}
[c]{ll}%
dY(t) & =-f(t,Y(t),E[Y(t)],Z(t))dt+Z(t)dB(t),t\in
\left[  0,T\right]  ,\\
Y(T) & =\xi,
\end{array}
\right.  \label{s-mfbsde}%
\end{equation}
where 
\begin{equation*}
f(t,Y(t),E[Y(t)],Z(t))=-r(t)Y(t)-r^{\prime}(t)E[Y(t)]+F(t,Z(t)).
\end{equation*}
\noindent We assume that the generator $(y, \bar{y}, z)\rightarrow f(t,Y(t),E[Y(t)],Z(t)):\Omega\times\lbrack0,T]\times\mathcal{
\mathbb{R}
}\times\mathcal{
\mathbb{R}
}\times\mathcal{
\mathbb{R}
}\rightarrow\mathcal{%
\mathbb{R}
}$ is $\mathbb{F}$-adapted, uniformly Lipschitz and concave, and the terminal condition
$\xi\in L^{2}\left(  \Omega,\mathcal{F}_{T}\right).$
\begin{definition}
Define $\varphi_{t}:$\emph{\ }%
$(T;\xi)\rightarrow\varphi_{t}(T;\xi)$ by%
\[%
\begin{array}
[c]{lll}%
\varphi_{t}(T;\xi) & = & -Y_{t}(T;\xi),\text{ \ \ \ \ }t\in\left[  0,T\right]
,
\end{array}
\]
where $Y_{t}(T;\xi)$ is a component of the solution of the MF-BSDE
(\ref{s-mfbsde}) with terminal horizon $T$, terminal condition $\xi$ and driver $f$. Then $\varphi_{t}(T;\xi)$ is a dynamic risk measure induced by a MF-BSDE.
\end{definition}

\noindent We may remark
that the driver $f$ depends linearly on $Y$ and its expected value $E[Y]$, and nonlinear with respect to $Z$. This is interpreted as a market with interest rates $(r(t),r'(t))$. We can reformulated this as a
problem with a driver independent of $Y$ and $E[Y]$ by discounting the financial position $\xi$. We assume that the instantaneous interest rates $r(t)$ and
$r^{\prime}(t)$ are deterministic. We denote by $\varphi_{\cdot}$, the
corresponding discounted risk-measure.\newline Define the discounted process%
\[
Y^{r}(t):=e^{-%
{\textstyle\int_{0}^{t}}
(r(s)+r^{\prime}(s))ds}Y(t).
\]
Then $Y^{r}$ with driver
\[
F^{r}(\cdot,t,Z(t)):=e^{-%
{\textstyle\int_{0}^{t}}
(r(s)+r^{\prime}(s))ds}F(\cdot,t,e^{-%
{\textstyle\int_{0}^{t}}
(r(s)+r^{\prime}(s))ds}Z(t)),
\]
and terminal value $\xi^{r}:=e^{-%
{\textstyle\int_{0}^{t}}
(r(s)+r^{\prime}(s))ds}\xi$ is a part of the solution of the associated BSDE.
We obtain also a discounted risk-measure accordingly
\[
\varphi_{0}(\xi,T)=\varphi_{0}^{r}(e^{-%
{\textstyle\int_{0}^{t}}
(r(s)+r^{\prime}(s))ds}\xi,T).
\]
This discounted risk-measure is translation-invariant because $F^{r}$ does not
depend on $Y$, we have for $\xi\in L^{2}\left(  \Omega,\mathcal{F}_{T}\right)
$ and $a\in%
\mathbb{R}
,$
\begin{align*}
\varphi_{0}(\xi+ae^{%
{\textstyle\int_{0}^{t}}
(r(s)+r^{\prime}(s))ds},T) &  =\varphi_{0}^{r}(e^{-%
{\textstyle\int_{0}^{t}}
(r(s)+r^{\prime}(s))ds}\xi+a,T)\\
&  =\varphi_{0}^{r}(e^{-%
{\textstyle\int_{0}^{t}}
(r(s)+r^{\prime}(s))ds}\xi,T)-a\\
&  =\varphi_{0}(\xi,T)-a.
\end{align*}
Similarly we can get for each $t\in\lbrack0,T]$, that
\[
\varphi(\xi,T)=\varphi^{r}(e^{-%
{\textstyle\int_{0}^{t}}
(r(s)+r^{\prime}(s))ds}\xi,T)
\]
is translation-invariant.

\subsection{Optimal portfolio with mean-field risk minimization}
Consider a financial market with two investment possibilities:\\
(i) Safe, or risk free asset with unit price 
\[
\
\begin{array}
[c]{ll}%
S_{0}(t)=1, \ t\in\lbrack0,T]\text{ .}\\
\end{array}
\]
(ii) Risky asset with unit price
\[
\
\begin{array}
[c]{ll}
dS_{1}(t)=S_{1}(t)[b_{0}(t)dt+\sigma_{0}(t)dB(t)], \ t\in\lbrack0,T]\text{ .}\\
\end{array}
\]
Let $\pi(t)$ be a self-financing portfolio invested in the risky asset at time $t$.
\noindent We want to minimize the risk $\varphi(X^{\pi}(T))$ of the terminal
value of the wealth process $X^{\pi}(t)$  corresponding to a portfolio $\pi$ which satisfies the linear SDE%

\begin{equation}
\left\{
\begin{array}
[c]{ll}%
dX^{\pi}(t) & =\pi(t)X^{\pi}(t)[b_{0}(t)dt+\sigma_{0}(t)dB(t)]\text{ , }%
t\in\lbrack0,T]\text{ ,}\\
X^{\pi}(0) & =x_{0}\text{ ,}%
\end{array}
\right.
\end{equation}

such that%
\[
\varphi(X^{\pi}(T))=-Y^{\pi}(0)
\]
where $Y^{\pi}(t)$ satisfies a MF-BSDE%
\begin{equation}
\left\{
\begin{array}
[c]{ll}%
-dY^{\pi}(t) & =[-r_{0}(t)E[Y^{\pi}(t)]+F(Z(t))]dt-Z(t)dB(t)\text{ , }%
t\in\lbrack0,T]\text{ ,}\\
Y^{\pi}(T) & =X^{\pi}(T)\text{ .}%
\end{array}
\right.  \label{MF_BSDE}%
\end{equation}
Here we assume that $b_{0}(t),$ $\sigma_{0}(t),$ $r_{0}(t)$ are given
deterministic functions and $F:%
\mathbb{R}
\rightarrow%
\mathbb{R}
$ is some given concave function. We want to find $\hat{\pi}\in\mathcal{U}%
_{\mathbb{G}}$ such that%
\[
\underset{\pi\in\mathcal{U}_{\mathbb{G}}}{\inf}\text{ }(-Y^{\pi}%
(0))=-Y^{\hat{\pi}}(0).
\]
Define the Hamiltonian $H$ that correspondds to our problem by%

\begin{align*}
H(t,x,z,\bar{y},\pi,p^{0},q^{0},\lambda^{0},\lambda^{1})  &  =p^{0}b_{0}\pi
x+q^{0}\sigma_{0}\pi x\\
&  +\lambda^{0}(r_{0}\bar{y}+F(z))+\langle\lambda^{1},\bar{y}\rangle.
\end{align*}
The couple $(p^{0},q^{0})$ solution of the following BSDE
\[
\left\{
\begin{array}
[c]{lll}%
dp^{0}(t) & = & -[p^{0}(t)b_{0}(t)\pi(t)+q^{0}(t)\sigma_{0}(t)\pi
(t)]dt+q^{0}(t)dB(t),t\in\left[  0,T\right]  ,\\
p^{0}(T) & = & \lambda^{0}(T),
\end{array}
\right.
\]
and $(p^{1},q^{1})$ satisfies%
\[
\left\{
\begin{array}
[c]{lll}%
dp^{1}(t) & = & q^{1}(t)dB(t),\ t\in\left[  0,T\right]  ,\\
p^{1}(T) & = & 0.
\end{array}
\right.
\]
$\lambda^{0}$ is given by the forward SDE%
\begin{equation}
\left\{
\begin{array}
[c]{lll}%
d\lambda^{0}(t) & = & \partial_{z}F(Z(t))\lambda^{0}(t)dB(t),\ t\in\left[
0,T\right]  ,\\
\lambda^{0}(0) & = & 1,
\end{array}
\right.  \label{lamda^0}%
\end{equation}
and
\[
\left\{
\begin{array}
[c]{lll}%
d\lambda^{1}(t) & = & -r_{0}(t)\lambda^{0}(t)dt,t\in\left[  0,T\right]  ,\\
\lambda^{1}(0) & = & 0.
\end{array}
\right.
\]
The first order necessary optimality condition gives%
\[
\hat{p}^{0}(t)b_{0}(t)\hat{X}(t)+\hat{q}^{0}(t)\sigma_{0}(t)\hat{X}(t)=0,
\]
where we denoted by $\hat{X}(t)=X^{\hat{\pi}}(t)$ and so on.\
Since $\hat{X}(t)>0$ for all $t$ $P$-a.s., we obtain
\begin{equation}
\hat{p}^{0}(t)b_{0}(t)+\hat{q}^{0}(t)\sigma_{0}(t)=0, \label{n_o_c}%
\end{equation}
which implies%
\[
\left\{
\begin{array}
[c]{lll}%
d\hat{p}^{0}(t) & = & \hat{q}^{0}(t)dB(t)=-\tfrac{b_{0}(t)}{\sigma_{0}(t)}%
\hat{p}^{0}(t)dB(t),t\in\left[  0,T\right]  ,\\
\hat{p}^{0}(T) & = & \hat{\lambda}^{0}(T),
\end{array}
\right.
\]
this together with equation (\ref{lamda^0}), yields%
\[
\hat{p}^{0}(t)=\hat{\lambda}^{0}(t)\text{, }\hat{q}^{0}(t)=\partial_{z}%
F(\hat{Z}(t))\hat{\lambda}^{0}(t).
\]
From (\ref{n_o_c}), we get%
\[
\partial_{z}F(\hat{Z}(t))=-\tfrac{b_{0}(t)}{\sigma_{0}(t)}.
\]
For example, if we choose%
\begin{equation}
F(z)=-\tfrac{1}{2}z^{2}. \label{cho}%
\end{equation}
That is
\[
\hat{Z}(t)=\tfrac{b_{0}(t)}{\sigma_{0}(t)}.
\]
Substituting the expression of $\hat{Z}(t)$ above into the MF-BSDE
(\ref{MF_BSDE}), we obtain%
\begin{equation}
\left\{
\begin{array}
[c]{ll}%
d\hat{Y}(t) & =-[-r_{0}(t)E[\hat{Y}(t)]-\tfrac{1}{2}(\tfrac{b_{0}(t)}%
{\sigma_{0}(t)})^{2}]dt-\tfrac{b_{0}(t)}{\sigma_{0}(t)}dB(t)\text{ , }%
t\in\lbrack0,T]\text{ ,}\\
\hat{Y}(T) & =\hat{X}(T)\text{ .}%
\end{array}
\right.  \label{y_^}%
\end{equation}
Consequently%
\[%
\begin{array}
[c]{ll}%
-dE[\hat{Y}(t)] & =[-r_{0}(t)E[\hat{Y}(t)]-\tfrac{1}{2}(\tfrac{b_{0}%
(t)}{\sigma_{0}(t)})^{2}]dt\text{ ,}%
\end{array}
\]
thus%
\begin{equation}%
\begin{array}
[c]{ll}%
E[\hat{Y}(t)] & =\exp(-%
{\textstyle\int_{0}^{t}}
r_{0}(s)ds)[\hat{Y}(0)+\tfrac{1}{2}%
{\textstyle\int_{0}^{t}}
\tfrac{b_{0}^{2}(s)}{\sigma_{0}^{2}(s)}\exp(%
{\textstyle\int_{0}^{s}}
r_{0}(\alpha)d\alpha)ds].
\end{array}
\label{y_bar}%
\end{equation}
Define $\Gamma(t)$ to be the solution of the linear SDE
\[
\left\{
\begin{array}
[c]{ll}%
d\Gamma(t) & =-\tfrac{b_{0}(t)}{\sigma_{0}(t)}\Gamma(t)dB(t)\text{ , }%
t\in\lbrack0,T]\text{ ,}\\
\Gamma(0) & =1\text{,}%
\end{array}
\right.
\]
or explicitely%
\begin{equation}
\Gamma(t)=\exp(-%
{\textstyle\int_{0}^{t}}
\tfrac{b_{0}(s)}{\sigma_{0}(s)}dB(s)-\tfrac{1}{2}%
{\textstyle\int_{0}^{t}}
(\tfrac{b_{0}(s)}{\sigma_{0}(s)})^{2}ds)\text{ , }t\in\lbrack0,T]\text{ .}
\label{ga}%
\end{equation}
By the Girsanov theorem of change of measures, we know that there exists an
equivalent local martingale measure $Q<<P$, such that%
\[
dQ=\Gamma(T)dP\text{ on }\mathcal{F}_{T},
\]
with $\Gamma(T)=\frac{dQ}{dP}$ is called the Radon-Nikodym derivative of $Q$
with respect to $P$ on $\mathcal{F}_{T}$.\newline Substituting (\ref{y_bar}%
)-(\ref{ga}) into (\ref{y_^}) we have \small
\begin{align*}
\hat{X}(T)  &  =\hat{Y}(T)=\hat{Y}(0)+\exp(-%
{\textstyle\int_{0}^{t}}
r_{0}(s)ds)[\hat{Y}(0)+\tfrac{1}{2}%
{\textstyle\int_{0}^{t}}
\tfrac{b_{0}^{2}(s)}{\sigma_{0}^{2}(s)}\exp(%
{\textstyle\int_{0}^{s}}
r_{0}(\alpha)d\alpha)]\\
&  +\tfrac{1}{2}%
{\textstyle\int_{0}^{T}}
(\tfrac{b_{0}(s)}{\sigma_{0}(s)})^{2}ds+%
{\textstyle\int_{0}^{T}}
\tfrac{b_{0}(s)}{\sigma_{0}(s)}dB(s)\\
&  =\hat{Y}(0)+\exp(-%
{\textstyle\int_{0}^{t}}
r_{0}(s)ds)[\hat{Y}(0)+\tfrac{1}{2}%
{\textstyle\int_{0}^{t}}
\tfrac{b_{0}^{2}(s)}{\sigma_{0}^{2}(s)}\exp(%
{\textstyle\int_{0}^{s}}
r_{0}(\alpha)d\alpha)ds]-\ln\Gamma(t).
\end{align*} \normalsize
Taking the expectation but now with respect to the new measure $Q$, we get%
\small
\begin{align}
-\hat{Y}(0)    =&-x_{0}-\exp(-%
{\textstyle\int_{0}^{t}}
r_{0}(s)ds)[\hat{Y}(0)+\tfrac{1}{2}%
{\textstyle\int_{0}^{t}}
\tfrac{b_{0}^{2}(s)}{\sigma_{0}^{2}(s)}\exp(%
{\textstyle\int_{0}^{s}}
r_{0}(\alpha)d\alpha)ds]-E_{Q}[\ln\Gamma(T)]\nonumber\\
  =&\tfrac{1}{1-\exp(-%
{\textstyle\int_{0}^{t}}
r_{0}(s)ds)}\{-x_{0}-\exp(-%
{\textstyle\int_{0}^{t}}
r_{0}(s)ds)[\tfrac{1}{2}%
{\textstyle\int_{0}^{t}}
\tfrac{b_{0}^{2}(s)}{\sigma_{0}^{2}(s)}\exp(%
{\textstyle\int_{0}^{s}}
r_{0}(\alpha)d\alpha)ds]\label{Y_hat}\\
  &- E[\Gamma(T)\ln\Gamma(T)],\nonumber
\end{align} \normalsize
where $E[\Gamma(T)\ln\Gamma(T)]$ is the entropy of $Q$ with respect to
$P$.\newline Since we obtained the optimal value of $\hat{Y}(0)$, we can get
the corresponding optimal terminal wealth $\hat{X}(T).$\newline Summarizing,
we have the following conclusion:

\begin{theorem}
Suppose that (\ref{cho}) holds. Then the minimal risk of our problem is given
by (\ref{Y_hat}).
\end{theorem}


\begin{bibdiv}
\begin{biblist}
\bib{andersson2011}{article}{
  title={A maximum principle for SDEs of mean-field type},
  author={Andersson, Daniel},
  author={Djehiche, Boualem},
  journal={Applied Mathematics \& Optimization},
  volume={63},
  number={3},
  pages={341--356},
  year={2011},
  publisher={Springer}
}


\bib{AO1}{article}{
title={Model uncertainty stochastic mean-field control},
  author={Agram, Nacira},
  author={{\O}ksendal, Bernt},
  journal={Stochastic Analysis and Applications},
  pages={1--21},
  year={2019},
  publisher={Taylor \& Francis}
}


\bib{bens}{article}{
  title={Well-posedness of mean-field type forward--backward stochastic differential equations},
  author={Bensoussan, Alain},
  author={Yam, Sheung Chi Phillip},
  author={Zhang, Zheng},
  journal={Stochastic Processes and their Applications},
  volume={125},
  number={9},
  pages={3327--3354},
  year={2015},
  publisher={Elsevier}
}


\bib{buckdahngeneral}{article}{
  title={A general stochastic maximum principle for SDEs of mean-field type},
  author={Buckdahn, Rainer},
  author={Djehiche, Boualem},
  author={Li, Juan},
  journal={Applied Mathematics \& Optimization},
  volume={64},
  number={2},
  pages={197--216},
  year={2011},
  publisher={Springer}
}



\bib{BDLP}{article}{
  title={Mean-field backward stochastic differential equations: a limit approach},
  author={Buckdahn, Rainer },
  author={Djehiche, Boualem},
  author={Li, Juan},
  author={Peng, Shige},
  journal={The Annals of Probability},
  volume={37},
  number={4},
  pages={1524--1565},
  year={2009},
  publisher={Institute of Mathematical Statistics}
}


\bib{carmona1}{article}{
  title={Control of McKean--Vlasov dynamics versus mean field games},
  author={Carmona, Ren{\'e}},
  author={Delarue, Fran{\c{c}}ois},
  author={Lachapelle, Aim{\'e}},
  journal={Mathematics and Financial Economics},
  pages={1--36},
  year={2013},
  publisher={Springer}
}


\bib{carmona2}{article}{
  title={Forward--backward stochastic differential equations and controlled McKean--Vlasov dynamics},
  author={Carmona, Ren{\'e}},
  author={Delarue, Fran{\c{c}}ois},
  journal={The Annals of Probability},
  volume={43},
  number={5},
  pages={2647--2700},
  year={2015},
  publisher={Institute of Mathematical Statistics}
}

\bib{Boualem2019}{article}{
  title={Mean-field backward-forward stochastic differential equations and nonzero sum stochastic differential games},
  author={Djehiche, Boualem},
  author={Hamadene, Said},
  journal={arXiv preprint arXiv:1904.06193},
  year={2019}
}


\bib{ll}{article}{
  title={Mean field games},
  author={Lasry, Jean-Michel},
  author={Lions, Pierre-Louis},
  journal={Japanese journal of mathematics},
  volume={2},
  number={1},
  pages={229--260},
  year={2007},
  publisher={Springer}
}

\bib{Lions}{misc}{
title={Cours au college de france: Th{\'e}orie des jeux \'a champs moyens},
  author={Lions, PL},
  year={2014}
}


\end{biblist}
\end{bibdiv}
\end{document}